\documentclass{article}
\usepackage[utf8]{inputenc}

\pdfoutput=1
\usepackage{mathtools}
\usepackage{amsfonts}
\usepackage{amssymb}
\usepackage{centernot}
\usepackage{graphicx}
\usepackage{amsthm}
\usepackage{enumerate}
\usepackage[noadjust]{cite}

\usepackage{tikz}
\usetikzlibrary{decorations.markings}
\usetikzlibrary{arrows}
\usetikzlibrary{calc}
\providecommand{\U}[1]{\protect\rule{.1in}{.1in}}
\newtheorem{theorem}{Theorem}

\newtheorem{corollary}[theorem]{Corollary}

\newtheorem{definition}[theorem]{Definition}

\newtheorem{lemma}[theorem]{Lemma}

\newtheorem{proposition}[theorem]{Proposition}

\newcommand{\Aut}{\textup{Aut}}

\newcommand{\As}{\textup{As}}
\newcommand{\tc}{\textup{tc}}

\theoremstyle{remark}

\begin{document}

\title{Linking numbers, quandles and groups}
\author{Lorenzo Traldi\\Lafayette College\\Easton, PA 18042, USA\\traldil@lafayette.edu
}
\date{ }
\maketitle

\begin{abstract}
We introduce a quandle invariant of classical and virtual links, denoted $Q _ {\tc} (L)$. This quandle has the property that $Q _ {\tc} (L) \cong Q _ {\tc} (L')$ if and only if the components of $L$ and $L'$ can be indexed in such a way that $L=K_1 \cup \dots \cup K_{\mu}$, $L'=K'_1 \cup \dots \cup K'_{\mu}$ and for each index $i$, there is a multiplier $\epsilon_i \in \{-1,1\}$ that connects virtual linking numbers over $K_i$ in $L$ to virtual linking numbers over $K'_i$ in $L'$: $\ell_{j/i}(K_i,K_j)= \epsilon_i \ell_{j/i}(K'_i,K'_j)$ for all $j \neq i$. We also extend to virtual links a classical theorem of Chen, which relates linking numbers to the nilpotent quotient $G(L)/G(L)_3$.


\end{abstract}

\section{Introduction}
Linking numbers are among the oldest, and simplest, invariants of classical knot theory; as discussed by Ricca and Nipoti \cite{RN}, Gauss mentioned linking numbers in his diary almost 200 years ago. Linking numbers were extended to virtual links by Goussarov, Polyak and Viro \cite{GPV}. For our purposes, the most convenient way to describe linking numbers is to calculate them from link diagrams.

In this paper the term \emph{link diagram} refers to an oriented virtual diagram; for a thorough discussion, we refer to the book of Manturov and Ilyutko \cite{MP}. Here is a summary. A link diagram is a subset $D$ of $\mathbb R ^2$ obtained from the union of a finite number $\mu$ of oriented, piecewise smooth closed curves. The curves must be in general position, i.e., there are only finitely many (self-) intersections among them, and all of these (self-) intersections are transverse double points (crossings). Each crossing is either \emph{classical} or \emph{virtual}; if no crossing is virtual, then $D$ is a \emph{classical} link diagram. The diagram is obtained from the union of the curves by: (a) at each classical crossing, removing a short piece of the underpassing segment on each side of the crossing, and (b) at each virtual crossing, drawing a small circle around the crossing. The result of removing the short pieces near classical crossings is to cut the original curves into arcs; the set of arcs in $D$ is denoted $A(D)$. These arcs are sometimes called the \emph{long arcs} of $D$, to contrast with shorter arcs obtained by cutting also at overpasses, or at virtual crossings. 

Two link diagrams are \emph{equivalent} or \emph{of the same type} if they are related through a finite sequence of four kinds of moves: detour moves affecting only virtual crossings, and classical Reidemeister moves affecting only classical crossings. Any of these moves provides a natural way to identify the closed curves underlying one diagram with the closed curves underlying the other diagram, so it makes sense to say that an equivalence class of link diagrams represents a $\mu$-component link type $L=K_1 \cup \dots \cup K_ \mu$. Reidemeister \cite{R} showed that equivalence classes of classical diagrams represents link types in $\mathbb S^3$, and Kuperberg \cite{K} showed that equivalence classes of virtual link diagrams represent link types in thickened surfaces.

\begin{figure} [bth]
\centering
\begin{tikzpicture} [>=angle 90]
\draw [thick] (1,0.5) -- (-0.6,-0.3);
\draw [thick] [<-] (-0.6,-0.3) -- (-1,-0.5);
\draw [thick] [->] (-1,0.5) -- (-.6,0.3);
\draw [thick] (-.6,0.3) -- (-.2,0.1);
\draw [thick] (0.2,-0.1) -- (1,-0.5);

\draw [thick] (6,0.5) -- (4.4,-0.3);
\draw [thick] [<-] (4.4,-0.3) -- (4,-0.5);
\draw [thick] (4,0.5) -- (4.4,0.3);
\draw [thick] [<-] (4.4,0.3) -- (4.8,0.1);
\draw [thick] (5.2,-0.1) -- (6,-0.5);
\node at (0,-1.4) {$w(c)=-1$};
\node at (5,-1.4) {$w(c)=1$};
\end{tikzpicture}
\caption{Classical crossings of negative and positive writhe.}
\label{writhefig}
\end{figure}
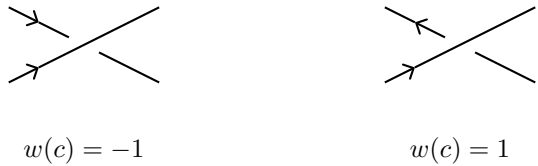

Linking numbers are defined using the notion of writhe of a classical crossing, illustrated in Fig. \ref{writhefig}.

\begin{definition}
\label{lno}
Let $D$ be a diagram of an oriented link, $L=K_1 \cup \dots \cup K_{\mu}$. If $i \neq j \in \{1, \dots , \mu\}$, let $C_{j/i}(D)$ be the set of classical crossings in $D$ at which $K_j$ passes over $K_i$. Then the \emph{linking number} of $K_j$ over $K_i$ is
\[
\ell_{j/i}(K_i,K_j) =\ell_{j/i}(K_j,K_i) = \sum_{c \in C_{j/i}(D)} w(c).
\]
\end{definition}

We refer to Ricca and Nipoti \cite{RN} for the history of the classical notion of linking numbers. The extension to virtual links given in Definition \ref{lno} was mentioned by Goussarov, Polyak and Viro \cite{GPV}. It is a simple exercise to verify that these linking numbers are invariant under Reidemeister moves and detour moves, so they are link type invariants. In the classical case, it is always true that $\ell_{i/j}(K_i,K_j) =\ell_{j/i}(K_i,K_j)$, so we can use the simpler notation $\ell(K_i,K_j) = \ell(K_j,K_i)$. For virtuals, instead, the two linking numbers $\ell_{i/j}(K_i,K_j)$ and $\ell_{j/i}(K_i,K_j)$ are independent of each other.

When Kauffman \cite{K} introduced virtual links, he extended many notions of classical knot theory to them. Among these notions were quandles, which had been introduced to classical knot theory by Joyce \cite{J} and Matveev \cite{M}.  We recall the definition.

\begin{definition}
\label{quandle}
A \emph{quandle} on a set $Q$ is specified by a binary operation $\triangleright$, with the following properties.
\begin{enumerate}
\item For each $x \in Q$, $x \triangleright x = x.$
\item For each $y \in Q$, the formula  $\beta_y(x)=x \triangleright y$ defines a bijection $\beta_y:Q \to Q$.
\item For all $w,x,y \in Q$, $(w \triangleright x) \triangleright y=(w \triangleright y) \triangleright (x \triangleright y)$.
\end{enumerate}
\end{definition}

If $Q$ and $Q'$ are quandles, then a function $f:Q \to Q'$  is a \emph{quandle map} if $f(x \triangleright y) = f(x) \triangleright f(y) \thickspace \allowbreak \forall x,y \in Q$. A bijective quandle map $Q \to Q'$ is an isomorphism, a bijective quandle map $Q \to Q$ is an automorphism, and the automorphisms of $Q$ form a group $\Aut(Q)$ under function composition. The $\beta_y$ maps mentioned in part 2 of Definition \ref{quandle} are the \emph{translations} of $Q$. Notice that the third requirement of Definition \ref{quandle} guarantees that the translations of $Q$ are automorphisms of $Q$. The subgroup of $\Aut(Q)$ generated by the translations is the \emph{translation group} of $Q$; we denote it $\beta(Q)$.

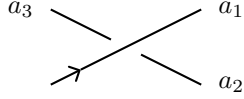
\begin{figure} [bth]
\centering
\begin{tikzpicture} [>=angle 90]
\draw [thick] (1,0.5) -- (-0.6,-0.3);
\draw [thick] [<-] (-0.6,-0.3) -- (-1,-0.5);
\draw [thick] (-1,0.5) -- (-.2,0.1);
\draw [thick] (0.2,-0.1) -- (1,-0.5);
\node at (1.4,0.5) {$a_1$};
\node at (-1.4,0.5) {$a_3$};
\node at (1.4,-0.5) {$a_2$};
\end{tikzpicture}
\caption{The arcs incident at a classical crossing.}
\label{crossfig}
\end{figure}

\begin{definition}
\label{linkquandle}
Let $D$ be a link diagram representing a link $L$. The \emph{quandle} of $L$ is the quandle generated by $\{q_a \mid a \in A(D)\}$, subject to the requirement that at each classical crossing as indicated in Fig.\ \ref{crossfig}, $q_{a_3} = q_{a_2} \triangleright q_{a_1}$. We denote this quandle $Q(L)$.
\end{definition}

It is not hard to see that $Q(L)$ is invariant under Reidemeister moves and detour moves, up to isomorphism, so it provides a link type invariant.

Here is a central definition of the paper.

\begin{definition}
\label{tcquandle}
A quandle $Q$ is \emph{translation-commutative} if its translation group $\beta(Q)$ is commutative. 
\end{definition}

As $\beta(Q)$ is generated by the translations, $Q$ is translation-commutative if and only if
\begin{equation}
\label{tc}
(x \triangleright y) \triangleright z = (x \triangleright z) \triangleright y \quad \thickspace \allowbreak \forall x, y, z \in Q.
\end{equation}
It follows that translation-commutative quandles constitute a variety of quandles, in the sense discussed by Joyce \cite[Sec.\ 10]{J}. Therefore, every quandle has a canonical translation-commutative quotient, described below.

\begin{definition}
\label{tcquotient}
Let $Q$ be a quandle. Then the \emph{translation-commutative quotient} of $Q$ is the quandle defined using generators and relations as follows. There is a generator $x_q$ for each $q \in Q$, there is a relation $x_{p \triangleright q} = x_p \triangleright x_q$ for each pair of elements $p,q \in Q$, and there is a relation $(x_p \triangleright x_q) \triangleright x_r = (x_p \triangleright x_r) \triangleright x_q$ for each triple of (not necessarily distinct) elements $p,q,r \in Q$. We denote this quandle $Q_{\tc}$.
\end{definition}

The mapping $q \mapsto x_q$ defines a canonical surjective quandle map $Q \to Q_{\tc}$. It is easy to see that if $S \subseteq Q$ generates $Q$, then the image of $S$ generates $Q_{\tc}$. 

In Sec.\ \ref{tcq}, we show that translation-commutative quandles are much less complicated than general quandles. A translation-commutative quandle is completely determined by a family of subgroups of a free abelian group, and it is not hard to determine whether or not two families of subgroups determine isomorphic translation-commutative quandles.

If $L$ is a link, then we define $Q_ {\tc}(L)$ to be the translation-commutative quotient of $Q(L)$, i.e., $Q_{\tc}(L) = Q(L)_{\tc}$. Our main result is that $Q_{\tc}(L)$ is strongly related to the linking numbers in $L$. 

\begin{theorem}
\label{thmone}
Suppose $L$ and $L'$ are links. Then $Q_{\tc}(L) \cong Q_{\tc}(L')$ if and only if both of the following requirements are satisfied. 
\begin{enumerate} [(a)]
    \item $L$ and $L'$ have the same number of components.
    \item The components of $L$ and $L'$ can be indexed so that $L=K_1 \cup \dots \cup K_\mu$, $L'=K'_1 \cup \dots \cup K'_\mu$, and there are $\epsilon_1, \dots, \epsilon_\mu \in \{-1,1\}$ such that $\ell_{j/i}(K_i,K_j) = \epsilon_i \ell_{j/i}(K'_i,K'_j)$ for all $i \neq j \in \{1, \dots, \mu \}$.
\end{enumerate}
\end{theorem}

Notice in particular that if $L$ and $L'$ are any knots (i.e., both have $\mu=1$), then the requirement regarding linking numbers in $L$ and $L'$ is satisfied vacuously, so Theorem \ref{thmone} asserts that $Q_{\tc}(L) \cong Q_{\tc}(L')$. In fact, the translation-commutative quandle of every knot is a trivial one-element quandle.

As mentioned above, every quandle has a canonical translation-commutative quotient. Therefore, Theorem \ref{thmone} directly implies the following.

\begin{corollary}
\label{med}
Let $\widetilde Q(L)$ be any invariant quandle of a link $L$, whose canonical translation-commutative quotient is $Q_{\tc}(L)$. If $L$ and $L'$ are links with $\widetilde Q(L) \cong \widetilde Q(L')$, then $L$ and $L'$ satisfy the two requirements of Theorem \ref{thmone}.
\end{corollary}

There are at least two well-known invariant quandles that can play the role of $\widetilde Q (L)$ in Corollary \ref{med}: the full link quandle, $Q(L)$, and the medial quandle, $\textup{AbQ}(L)$ in Joyce's notation \cite{J}.

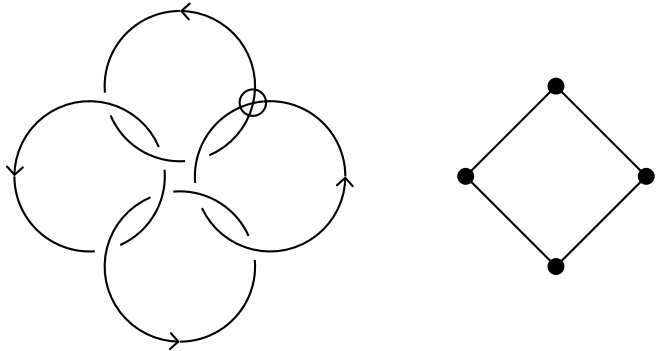
\begin{figure} [bth]
\centering
\begin{tikzpicture} [>=angle 90]
\draw [->] [thick, domain=-67:90] plot ({cos(\x)}, {sin(\x)});
\draw [thick, domain=90:185] plot ({cos(\x)}, {sin(\x)});
\draw [thick, domain=-86:-157]  plot ({cos(\x)}, {sin(\x)});

\draw [thick, domain=-66:5] plot ({-1.2+cos(\x)}, {-1.2+sin(\x)});
\draw [->] [thick, domain=23:180] plot ({-1.2+cos(\x)}, {-1.2+sin(\x)});
\draw [thick, domain=180:274] plot ({-1.2+cos(\x)}, {-1.2+sin(\x)});

\draw [thick, domain=114:185] plot ({1.2+cos(\x)}, {-1.2+sin(\x)});
\draw [->] [thick, domain=-155:0] plot ({1.2+cos(\x)}, {-1.2+sin(\x)});
\draw [thick, domain=0:95] plot ({1.2+cos(\x)}, {-1.2+sin(\x)});

\draw [thick, domain=114:95] plot ({1.2+cos(\x)}, {-1.2+sin(\x)});
\draw [thick] (.97,-.22) circle (5pt);

\draw [->] [thick, domain=-247:-90] plot ({cos(\x)}, {-2.4+sin(\x)});
\draw [thick, domain=-90:5] plot ({cos(\x)}, {-2.4+sin(\x)});
\draw [thick, domain=24:95]  plot ({cos(\x)}, {-2.4+sin(\x)});
\filldraw [black] (5,0) circle (3pt);
\filldraw [black] (3.8,-1.2) circle (3pt);
\filldraw [black] (6.2,-1.2) circle (3pt);
\filldraw [black] (5,-2.4) circle (3pt);
\draw [thick] (5,0) -- (3.8,-1.2);
\draw [thick] (5,0) -- (6.2,-1.2);
\draw [thick] (5,-2.4) -- (3.8,-1.2);
\draw [thick] (5,-2.4) -- (6.2,-1.2);
\end{tikzpicture}
\caption{A four-component link, and its linking graph.}
\label{insepfig}
\end{figure}

In general, the sign changes $\epsilon_1, \dots, \epsilon_\mu$ in Theorem \ref{thmone} can vary independently. When $L$ and $L'$ are both classical links, however, the sign changes are tied together: if $\ell(K_i,K_j) \neq 0$, then $\epsilon_i$ must equal $\epsilon_j$. One way to make this interdependence explicit involves the following construction. See Fig.\ \ref{insepfig} for an example. 

\begin{definition}
\label{lg}
A link $L=K_1 \cup \dots \cup K_{\mu}$ has a \emph{linking graph} $ \ell g(L)$, defined as follows.
\begin{itemize}
    \item There is a vertex $v_i$ for each component $K_i$ of $L$ such that $\ell_{j/i}(K_i,K_j) \neq 0$ or $\ell_{i/j}(K_i,K_j) \neq 0$ for some $j \neq i \in \{1, \dots, \mu \}$.
    \item Two vertices $v_i$ and $v_j$ are adjacent if $i \neq j$ and $\ell_{j/i}(K_i,K_j) \neq 0$ or $\ell_{i/j}(K_i,K_j) \neq 0$.
\end{itemize}
\end{definition}

Note that in general, the number of vertices in $ \ell g(L)$ may be less than $\mu$. (Indeed, $\ell g(L)$ may be empty.) Also, $ \ell g(L)$ is a simple graph; there are neither loops nor parallel edges.

Recall that a \emph{connected component} of a graph is an equivalence class of vertices under the equivalence relation $\sim$ generated by $v_i \sim v_j$ whenever $v_i$ and $ v_j$ are adjacent. We always use the phrase ``connected component'' for this notion, to avoid confusion with the components of a link. A graph with only one connected component is \emph{connected}. We consider the empty graph to be connected.

As classical links always have $\ell_{j/i}(K_i,K_j)=\ell_{i/j}(K_i,K_j)$, Theorem \ref{thmone} immediately implies the following.

\begin{theorem}
\label{thmtwo}
Let $L$ and $L'$ be classical links. Then $Q_{\tc}(L) \cong Q_{\tc}(L')$ if and only if both of the following conditions hold.
\begin{enumerate} [(a)]
    \item $L$ and $L'$ have the same number of components.
    \item The components of $L$ and $L'$ can be indexed in such a way that for each connected component $C$ of $\ell g(L)$, there is an $\epsilon_C \in \{-1,1\}$ such that $\ell(K_i,K_j) = \epsilon_C\ell(K'_i,K'_j)$ whenever $v_i$ is a vertex of $C$.
\end{enumerate}
\end{theorem}

In particular, if $\ell g (L)$ has only one connected component then (b) requires an $\epsilon \in \{-1,1\}$ such that $\ell(K_i,K_j) = \epsilon\ell(K'_i,K'_j) \thickspace \allowbreak \forall i \neq j \in \{1, \dots, \mu \}$.

In the classical case, Theorems \ref{thmone} and \ref{thmtwo} are the strongest results one could hope for. Even the full quandle $Q(L)$ of a classical link -- a very sensitive link invariant -- does not detect linking numbers absolutely. The reason is that $Q(L)$ is invariant under the combination of reversing the orientations of all the components of $L$, and replacing $L$ with its mirror image. Moreover, for a split link $L$, $Q(L)$ is not changed if the component orientations are reversed in just one split portion of $L$, and that portion is replaced with its mirror image. The effect of replacing a split portion with the orientation-reversed mirror image is to multiply all linking numbers from the connected component(s) of $\ell g (L)$ corresponding to that portion of the link by $-1$. 

Here is an outline of the paper. In Sec.\ \ref{tcq} we provide a complete structure theory for translation-commutative quandles, and in Sec.\ \ref{qtc} we use this structure theory to prove Theorem \ref{thmone}. In Sec.\ \ref{clab} we relate $Q_{\tc}(L)$ to the nilpotent quotient $G(L)/G(L)_3$ of the group of $L$. In Sec.\ \ref{hn} we briefly discuss the connection between our results and those of Harrell and Nelson \cite{HN}, who related linking numbers of two-component links to quandle counting invariants.

\section{Translation-commutative quandles}
\label{tcq}

In this section, we give a structure theory for translation-commutative quandles. If $B$ is a nonempty set, we use $\mathbb Z _B$ to denote the free abelian group on $B$.

\begin{definition}
\label{qs}
Let $S=\{S_b \mid b \in B\}$ be a family of subgroups of $\mathbb Z _B$, with $b \in S_b \thickspace \allowbreak \forall b \in B$. For each $b \in B$, let $A_b = \mathbb Z _B / S_b$. Then $Q (S)$ is the disjoint union 
\[
Q(S) = \bigcup \limits _{b \in B} A_b \text{,}
\]
equipped with the operation $\triangleright$ defined as follows: if $x \in A_b$ and $y \in A_c$, then 
\[
x \triangleright y = x+(c+S_b) \in A_b.
\]
\end{definition}

As a minor abuse of notation we often write $x \triangleright y =x+c+S_b$, even though $x+c$ is not well defined.

\begin{proposition}
\label{fgtc} $Q (S)$ is a translation-commutative quandle.
\end{proposition}
\begin{proof}
If $x \in A_b$ then as $b \in S_b$,
\[
x \triangleright x = x+b +S_b = x + 0 + S_b =x.
\]

If $y \in A_c$ then for each $b \in B$, the map $A_b \to A_b$ defined by $x \mapsto x \triangleright y$ is a bijection, with inverse function given by $x \mapsto x-(c+S_b) \in A_b$.

If $x \in A_b$, $y \in A_c$ and $z \in A_d$, then $y \triangleright z \in A_c$, so
\[
(x \triangleright y) \triangleright z = x+c+d + S_b = x+d+c+S_b
\]
\[
= (x+d+S_b) \triangleright (y \triangleright z) =  (x \triangleright z) \triangleright (y \triangleright z).
\]
Also, formula (\ref{tc}) is satisfied because
\[
(x \triangleright y) \triangleright z= x+c+d + S_b = x+d+c+S_b
= (x \triangleright z) \triangleright y.
\] 
\end{proof}

Recall that if $x$ is an element of a quandle $Q$, then the \emph{orbit} of $x$ in $Q$ is the smallest subset $Q_x \subseteq Q$ such that $x \in Q_x$ and $\beta_y^{\pm 1}(z) \in Q_x \thickspace \allowbreak \forall y \in Q \thickspace \allowbreak \forall z \in Q_x$. It is easy to see that for a quandle of type $Q (S)$, the orbits are the subsets $A_b$. Also, Definition \ref{qs} implies that if $b \in B$ and $y,z \in A_b$, then $\beta_y = \beta_z$. We use $\beta_b$ to denote the map $\beta_y$ for all $y \in A_b$.

\begin{lemma}
\label{sker}
Suppose $n \in \mathbb N$, $b_1, \dots, b_n \in B$, and $m_1, \dots, m_n \in \mathbb Z$. Let
\[
y = \sum_{i=1}^n m_i b_i \in \mathbb Z_B \text{,}
\]
and let 
\[
\beta_y = \prod_{i=1}^n \beta ^{m_i}_ {b_i} \in \beta(Q(S)).
\]
Then for each $b \in B$, the following are equivalent to each other.
\begin{enumerate}
    \item $y \in S_b$.
    \item For some $x \in A_b$, $\beta_y(x) = x$.
    \item For all $x \in A_b$, $\beta_y(x) = x$.
\end{enumerate}
\end{lemma}
\begin{proof}
According to Definition \ref{qs}, if $x \in A_b$ then
\[
\beta_y(x) = \beta_{b_1}^{m_1} \cdots \beta_{b_n}^{m_n}(x) = x+y+S_b.
\]
Therefore, $\beta_y(x) = x$ if and only if $y+S_b=S_b$.
\end{proof}
If $f:B \to B'$ is a function of sets, then we also use $f$ to denote the homomorphism $\mathbb Z _B \to \mathbb Z _{B'}$ defined by $f$. 
\begin{proposition}
\label{iso1}
Suppose $B$ and $B'$ are nonempty sets, $\{S_b \mid b \in B\}$ is a family of subgroups of $\mathbb Z _B$ with $b \in S_b \thickspace \allowbreak \forall b \in B$, and $\{S'_{b'} \mid b' \in B'\}$ is a family of subgroups of $\mathbb Z _{B'}$ with $b' \in S'_{b'} \thickspace \allowbreak \forall b' \in B'$. Then $Q(S) \cong Q(S')$ if and only if there is a bijection $f:B \to B'$ with $f(S_b)=S'_{f(b)} \thickspace \allowbreak \forall b \in B$.
\end{proposition}
\begin{proof}
We modify notation used for $Q(S)$ by using apostrophes when discussing $Q(S')$. For instance, if $b' \in B'$ then $A'_{b'}  = \mathbb Z _{B'} / S'_{b'}$. If there is a bijection $f$ as described, then for each $b \in B$, $f$ defines a group isomorphism $A_b \to A'_{f(b)}$. These group isomorphisms define a quandle isomorphism $Q(S) \to Q(S')$.

For the converse, suppose $g:Q(S) \to Q(S')$ is a quandle isomorphism. As an isomorphism, $g$ must define a bijection between the orbits of $Q (S)$ and the orbits of $Q (S')$, so $g$ must define a bijection $f:B \to B'$ with $g(A_b) = A'_{f(b)} \thickspace \allowbreak \forall b \in B$. 

As $g$ is a quandle isomorphism, $g(x \triangleright y) = g(x) \triangleright g(y) \thickspace \allowbreak \forall x,y \in Q(S)$. It follows that $g \circ \beta_b ^{\pm 1} = (\beta'_{f(b)})^{\pm 1} \circ g \thickspace \allowbreak \forall b \in B$, and hence 
\begin{equation}
    \label{comm}
g \circ \beta_{b_1}^{m_1} \cdots \beta_{b_n}^{m_n} = (\beta'_{f(b_1)})^{m_1} \cdots (\beta'_{f(b_n)})^{m_n} \circ g 
\end{equation}
whenever $n \in \mathbb N$, $ b_1, \dots, b_n \in B$ and $m_1, \dots, m_n \in \mathbb Z$.
Now, suppose $b \in B$. According to (\ref{comm}), if $x \in A_b$, $ b_1, \dots, b_n \in B$ and $m_1, \dots, m_n \in \mathbb Z$ then $\beta_{b_1}^{m_1} \cdots \beta_{b_n}^{m_n}(x)=x$ if and only if $(\beta'_{f(b_1)})^{m_1} \cdots (\beta'_{f(b_n)})^{m_n}(g(x)) = g(x)$. We deduce from Lemma \ref{sker} that $\sum m_ib_i \in S_b$ if and only if $\sum m_i f(b_i) \in S'_{f(b)}$. That is, $f(S_b) = S'_{f(b)}$. 
\end{proof}

A special case of Proposition \ref{iso1} will be particularly useful. Observe that if $C \subset B$ then there is a natural way to consider $\mathbb Z _C$ as a subgroup of $\mathbb Z _B$: simply identify each element $\sum n_i c_i \in \mathbb Z_C$ with the element $\sum n_i c_i \in \mathbb Z_B$ given by the same linear combination of generators..

\begin{corollary}
\label{iso}
Suppose $B$ and $B'$ are nonempty sets, $s_b \in \mathbb Z_{B-b}  \thickspace \allowbreak \forall b \in B$, and $s'_{b'} \in \mathbb Z_{B'-b'}  \thickspace \allowbreak \forall b' \in B'$. For each $b \in B$, let $S_b$ be the subgroup of $\mathbb Z _B$ generated by $\{b,s_b\}$, and for each $b' \in B'$, let $S'_{b'}$ be the subgroup of $\mathbb Z _{B'}$ generated by $\{b',s'_{b'}\}$. Then $Q(S) \cong Q(S')$ if and only if there is a bijection $f:B \to B'$ with $f(s_b)= \pm s'_{f(b)} \thickspace \allowbreak \forall b \in B$.
\end{corollary}
\begin{proof}
If $b \in B$ and $x,y \in \mathbb Z_{B-b}$, then $\{b,x\}$ and $\{b,y\}$ generate the same subgroup of $\mathbb Z_B$ if and only if $x = \pm y$. The corollary follows from this property and Proposition \ref{iso1}.
\end{proof}
The next result will help us show that every translation-commutative quandle is isomorphic to a quandle of type $Q(S)$.

\begin{lemma}
\label{betamaps}
Suppose $Q$ is a translation-commutative quandle, and $y \in Q$. Then $\beta_y = \beta_z \thickspace \allowbreak \forall z \in Q_y$.
\end{lemma}
\begin{proof}
If $y \in Q$ then according to Definition \ref{quandle}, for any $w,x \in Q$
\[
\beta_{\beta_x(y)}(w) = w \triangleright (y \triangleright x) = (\beta_x^{-1}(w) \triangleright x) \triangleright (y \triangleright x) = (\beta_x^{-1}(w) \triangleright y)  \triangleright x.
\]
According to formula (\ref{tc}), it follows that for every $w \in Q$,
\[
\beta_{\beta_x(y)}(w) =  (\beta_x^{-1}(w) \triangleright x)  \triangleright y = w \triangleright y = \beta_y(w).
\]
Hence $\beta_{\beta_x(y)}=\beta_y$. \end{proof}

\begin{theorem}
\label{tcstruc}
Let $Q$ be a translation-commutative quandle, and suppose there is a bijection $b \leftrightarrow Q_b$ between a set $B$ and the set of orbits of $Q$. Then there is a family $S = \{S_b \mid b \in B \}$ of  subgroups of $\mathbb Z _B$ such that $Q \cong Q(S)$.
\end{theorem}
\begin{proof}
For each $b \in B$, choose a representative element $q_b \in Q_b$, and let $\beta_b = \beta_{q_b}$. According to Lemma \ref{betamaps}, $\{ \beta_b \mid b \in B \}$ includes all the $\beta$ maps of $Q$, so $\beta(Q)$ is an abelian group generated by $\{ \beta_b \mid b \in B \}$. For each $b \in B$, then, there is a surjective function $f_b:\mathbb Z _B \to Q_b$ given by 
\[
f_b \left(\sum_{i=1}^n m_i b_i \right)= \left( \prod_{i=1}^n \beta_{b_i}^{m_i} \right) (q_b)
\]
whenever $n \geq 0 \in \mathbb Z$, $ b_1, \dots, b_n \in B$ and $m_1, \dots, m_n \in \mathbb Z$.

For each $b \in B$, let $S_b = \{\sum m_i b_i \in \mathbb Z _B \mid f_b(\sum m_i b_i) = q_b\}$. The fact that $Q$ is translation-commutative implies that $S_b$ is a subgroup of $\mathbb Z _B$, and $f_b$ induces a bijection between $\mathbb Z _B /S_b$ and $Q_b$. Also, the fact that $\beta_b(q_b) = \beta_{q_b}(q_b) = q_b$ implies that $b \in S_b$. It follows that taken together, the surjections $f_b$ induce an isomorphism $Q(S) \to Q$.
\end{proof}

\section{Theorems \ref{thmone} and \ref{thmtwo}}
\label{qtc}

Let $D$ be a diagram of an oriented link $L=K_1 \cup \dots \cup K_\mu$. The image in $D$ of each component $K_i$ consists of finitely many arcs, separated at classical crossings where $K_i$ is the underpassing component. For each $i$, choose a fixed arc $a_{i0}$ belonging to the image of $K_i$ in $D$.

If $a \in A(D)$ then for convenience, we use $x_a$ to denote the element $x_{q_a} \in Q_{\tc}(L)$.

\begin{lemma}
\label{orb}
The quandle $Q_{\tc}(L)$ has $\mu$ orbits. For each $i \in \{1, \dots, \mu\}$, there is an orbit that contains all the elements $x_a$ such that $a \in A(D)$ belongs to the image of $K_i$ in $D$. This orbit does not contain any element $x_a$ such that $a \in A(D)$ belong to the image of $K_j$ in $D$, where $j \neq i$.
\end{lemma}
\begin{proof} Consider a crossing of $D$ as pictured in Fig.\ \ref{crossfig}. 
Definitions \ref{qtc} and \ref{tcquotient} tell us that $x_{a_3} = x_{a_2} \triangleright x_{a_1}$, so $x_{a_2}$ and $x_{a_3}$ belong to the same orbit of $Q_{\tc}(L)$. 

Now, choose any $i \in \{1, \dots, \mu\}$. If we walk along the image of $K_i$ in $D$, starting at $a_{i0}$, then we ultimately encounter every arc of $D$ belonging to $K_i$. The observation of the previous paragraph applies each time we pass from one arc of $D$ to a different arc of $D$, so the $x_a$ elements corresponding to arcs of $K_i$ all belong to the same orbit of $Q_{\tc}(L)$.

To show that these orbits are distinct,  consider the trivial quotient quandle $Q_0$ of $Q_{\tc}(L)$, i.e., the quandle obtained by adding $x_p \triangleright x_q = x_p \allowbreak \thickspace \forall p,q \in Q$ to the requirements of Definition \ref{tcquotient}, with $Q=Q(L)$. It is clear that the definition of this quotient quandle is satisfied by the set $\{1, \dots, \mu\}$, considered as a trivial quandle, with all $\beta$ maps equal to the identity map. The fact that there is a well-defined quotient map $Q_{\tc}(L) \to Q_0$ guarantees that $Q_{\tc}(L)$ has $\mu$ different orbits.
\end{proof}

For each $i \in \{1, \dots, \mu\}$, let $\beta_i = \beta_{x_{a_{i0}}}:Q_{\tc}(L) \to Q_{\tc}(L)$. According to Lemma \ref{betamaps}, $\beta_1, \dots, \beta_\mu$ are all the $\beta$ maps of $Q_{\tc}(L)$. It follows that $Q_{\tc}(L)$ can be described like this: 

\begin{proposition}
\label{qtc2}$Q_{\tc}(L)$ is the translation-commutative quandle generated by $\{x_a \mid a \in A(D)\}$, subject to these requirements.
\begin{itemize}
    \item $Q_{\tc}(L)$ has $\beta$ maps $\beta_1, \dots, \beta_\mu$, corresponding to the components of $L$.
    \item At each classical crossing as indicated in Fig.\ \ref{crossfig}, if the overpassing arc $a_1$ belongs to the image of $K_\kappa$ in $D$, then $x_{a_3} = \beta_\kappa(x_{a_2})$.
\end{itemize} 
\end{proposition}

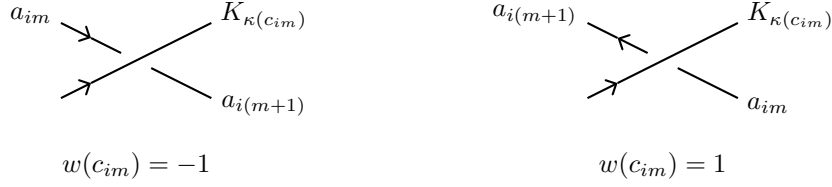
\begin{figure} [bth]
\centering
\begin{tikzpicture} [>=angle 90]
\draw [thick] (1,0.5) -- (-0.6,-0.3);
\draw [thick] [<-] (-0.6,-0.3) -- (-1,-0.5);
\draw [thick] [->] (-1,0.5) -- (-.6,0.3);
\draw [thick] (-.6,0.3) -- (-.2,0.1);
\draw [thick] (0.2,-0.1) -- (1,-0.5);

\draw [thick] (8,0.5) -- (6.4,-0.3);
\draw [thick] [<-] (6.4,-0.3) -- (6,-0.5);
\draw [thick] (6,0.5) -- (6.4,0.3);
\draw [thick] [<-] (6.4,0.3) -- (6.8,0.1);
\draw [thick] (7.2,-0.1) -- (8,-0.5);
\node at (-1.4,.6) {$a_{im}$};
\node at (1.7,-0.6) {$a_{i(m+1)}$};
\node at (1.7,0.6) {$K_{\kappa(c_{im})}$};
\node at (8.4,-0.6) {$a_{im}$};
\node at (8.7,0.6) {$K_{\kappa(c_{im})}$};
\node at (5.3,.6) {$a_{i(m+1)}$};
\node at (0,-1.4) {$w(c_{im})=-1$};
\node at (7,-1.4) {$w(c_{im})=1$};
\end{tikzpicture}
\caption{$K_i$ passes under $K_{\kappa(c_{im})}$ at the classical crossing $c_{im}$.}
\label{writhefig2}
\end{figure}

To describe $Q_{\tc}(L)$ in greater detail, it is convenient to index the arcs and classical crossings of $D$. Suppose that for each $i \in \{1, \dots, \mu \}$, $n_i$ is the number of arcs of $D$ belonging to $K_i$. We index these arcs $a_{i0}, a_{i1}, \dots, a_{i(n_i-1)}, a_{in_i} = a_{i0}$ in order, as we walk along $K_i$ in the direction of its orientation. If $0 \leq m < n_i$, let $c_{im}$ be the classical crossing we traverse as we walk from $a_{im}$ to $a_{i(m+1)}$, let $K_{\kappa(c_{im})}$ be the overpassing component of $L$ at $c_{im}$, and let $w(c_{im}) \in \{-1,1\}$ be the writhe of $c_{im}$. See Fig.\ \ref{writhefig2}.

Proposition \ref{qtc2} implies that if $w(c_{im})=-1$ then $x_{a_{im}} = \beta_ {\kappa(c_{im})}(x_{a_{i(m+1)}})$, and if $w(c_{im})=1$ then $x_{a_{i(m+1)}} = \beta_ {\kappa(c_{im})}(x_{a_{im}})$. In either case, $x_{a_{i(m+1)}} = \beta_ {\kappa(c_{im})}^{w(c_{im})}(x_{a_{im}})$. Therefore if $1 \leq i \leq \mu$ and $0\leq m \leq n_i-1$,
\begin{equation}
\label{crosses}
x_{a_{i(m+1)}} = \left( \prod _{k=0}^{m}\beta_ {\kappa(c_{ik})} ^{w(c_{ik})} \right) (x_{a_{i0}}).
\end{equation}
An instance of (\ref{crosses}) with $0 \leq m\leq n_i-2$ allows us to express the quandle element $x_{a_{i(m+1)}}$ in terms of $x_{a_{10}}, \dots, x_{a_{ \mu 0}}$. After we remove all of these redundant $x_{a_{i(m+1)}}$ generators, only $x_{a_{10}}, \dots, x_{a_{ \mu 0}}$ remain. We rename these elements $x_1, \dots, x_ \mu$, for simplicity. 

The instances of (\ref{crosses}) that have not been used to remove redundant generators are the ones with $m=n_i-1$. They tell us that for $1 \leq i \leq \mu$,
\begin{equation}
\label{defform}
x_i = \left( \prod _{\substack{j=1 \\ j \neq i}}^{\mu}\beta_ {j} ^{\ell_{j/i}(K_i,K_{j})} \right) (x_i).
\end{equation}

We deduce yet another equivalent description of $Q_{\tc}(L)$:
\begin{proposition}
\label{qtc3}
$Q_{\tc}(L)$ is the translation-commutative quandle generated by $\{x_1, \dots, x_ \mu\}$, subject to the requirement that (\ref{defform}) holds for every $i \in \{1, \dots, \mu \}$.
\end{proposition}

Let $B=\{b_1, \dots, b_\mu \}$, and for each $i \in \{1, \dots, \mu\}$ let $S_i(L)=S_{b_i}(L)$ be the subgroup of $\mathbb Z _B$ generated by $b_i$ and $\ell_i$, where
\begin{equation}
\label{long}
\ell_i = \sum _{\substack{j=1 \\ j \neq i}}^{\mu} \ell_{j/i}(K_i,K_{j}) b_{j}.
\end{equation}
Let $S(L)=\{S_{b_i}(L) \mid  1\leq i \leq \mu \}$. Then Proposition \ref{qtc3} implies the following.
\begin{corollary}
\label{qtc4}
The quandles $Q_{\tc}(L)$ and $Q(S(L))$ are isomorphic.
\end{corollary}


Theorem \ref{thmone} follows from Corollaries \ref{iso} and \ref{qtc4}. As discussed in the introduction, Theorem \ref{thmtwo} follows directly from Theorem \ref{thmone} and the fact that for a classical link, $\ell_{j/i}(K_i,K_{j})  = \ell_{i/j}(K_i,K_{j}) \thickspace \allowbreak \forall i \neq j \in \{1, \dots, \mu \}$.

\section{The quotient $G(L)/G(L)_3$}
\label{clab}

In this section we connect the ideas of the preceding sections to a theorem of Chen \cite{C} involving the nilpotent quotient $G(L)/G(L)_3$, where $G(L)= \pi_1(\mathbb S^3 - L)$ is the fundamental group of the complement of a classical link. We begin with a construction that dates back to the introduction of quandles to knot theory by Joyce \cite{J} and Matveev \cite{M}.

\begin{definition}
\label{As}
If $Q$ is a quandle, then there is an associated group $\As(Q)$, described using generators and relations as follows. There is a generator $g_ q$ for each $q \in Q$, and for each pair $p,q \in Q$ there is a relation $g_{p \triangleright q} = g_ q g_ p g_ q ^{-1}$.
\end{definition}

We should mention that Joyce \cite{J} used the notation $\textup{Adconj}(Q)$ instead of $\As(Q)$. Also, his quandle is the ``opposite'' of ours, i.e., it is defined by $g_{p \triangleright q} = g_ q ^{-1} g_ p g_ q $ instead of $g_{p \triangleright q} = g_ q g_ p g_ q ^{-1}$. We use the convention of Definition \ref{As} so that if $L$ is a link with a diagram $D$, then the presentation of $\As(Q(L))$ provided by Definitions \ref{linkquandle} and \ref{As} is the same as the presentation of the group $G(L)$ given by Kauffman \cite{K}. In the classical case, the same presentation was given by Fox \cite{F}. 

Recall that if $G$ is a group and $H$ is a subgroup of $G$, then $[G,H]$ is the subgroup of $G$ generated by the set of commutators $[g,h]=ghg^{-1}h^{-1}$ with $g \in G$ and $h \in H$. The \emph{lower central series} of $G$ is the sequence of normal subgroups $G=G_1 \supseteq G_2 \supseteq \dots $ with $G_{n+1} = [G,G_n]$ for each $n \geq 2$, and $G$ is \emph{nilpotent of class} $n$ if $G_n \neq \{1\} = G_{n+1}$. In particular, nontrivial abelian groups are nilpotent of class 1.

\begin{proposition}
\label{nilpotent}
Let $Q$ be a quandle. Then the quandle surjection $Q \to Q_{\tc}$ induces a group surjection $\As(Q) \to \As(Q_{\tc})$, whose kernel is the lower central series subgroup $\As(Q)_3$.
\end{proposition}
\begin{proof}
The kernel of the induced surjection $\As(Q) \to \As(Q_{\tc})$ is the normal subgroup of $\As(Q)$ generated by the images of relators derived from the quandle relations $(x_p \triangleright x_q) \triangleright x_r = (x_p \triangleright x_r) \triangleright x_q$ for $p,q,r \in Q$.

The quandle relation $(x_p \triangleright x_q) \triangleright x_r = (x_p \triangleright x_r) \triangleright x_q$ provides the group relation $g_{x_r}g_{x_q}g_{x_p}g_{x_q}^{-1}g_{x_r}^{-1}=g_{x_q}g_{x_r}g_{x_p}g_{x_r}^{-1}g_{x_q}^{-1}$, which corresponds to the relator
\[
g_{x_r}g_{x_q}g_{x_p}g_{x_q}^{-1}g_{x_r}^{-1}g_{x_q}g_{x_r}g_{x_p}^{-1}g_{x_r}^{-1}g_{x_q}^{-1}.
\]
Conjugating by $g_{x_p}^{-1}g_{x_r}^{-1}g_{x_q}^{-1}$, we see that this relator is equivalent to
\[
g_{x_p}^{-1}g_{x_r}^{-1}g_{x_q}^{-1}g_{x_r}g_{x_q}g_{x_p}g_{x_q}^{-1}g_{x_r}^{-1}g_{x_q}g_{x_r}=[g_{x_p}^{-1},[g_{x_r}^{-1},g_{x_q}^{-1}]].
\]
As $\As(Q)$ is generated by the inverses of the elements $g_{x_y}$ for $y \in Q$, the commutators $[g_{x_p}^{-1},[g_{x_r}^{-1},g_{x_q}^{-1}]]$ generate $\As(Q)_3$ as a normal subgroup of $\As(Q)$.
\end{proof}

\begin{corollary}
\label{linknilpotent}
Let $L$ be a link with group $G(L)$. Then the nilpotent quotient group $G(L)/G(L)_3$ has a  presentation with a generator $g_i$ for each $i \in \{1, \dots, \mu \}$ and two kinds of relations. 
\begin{enumerate}
    \item If $i \in \{1, \dots, \mu \}$, there is a relation
    \[
g_i = \left( \prod _{\substack{j=1 \\ j \neq i}}^{\mu}g_ {j} ^{\ell_{j/i}(K_i,K_{j})} \right)  g_i  \left( \prod _{\substack{j=1 \\ j \neq i}}^{\mu}g_ {j} ^{\ell_{j/i}(K_i,K_{j})} \right)^{-1}.
\]
    \item If $i,j,k \in \{1, \dots,\mu\}$,  there is a relation $[g_i,[g_j,g_k]]=1$.
    \end{enumerate} 

\end{corollary}
\begin{proof}
Proposition \ref{nilpotent} implies that the nilpotent quotient $G(L)/G(L)_3$ is isomorphic to $\As(Q_{\tc}(L))$. According to Proposition \ref{qtc3}, the quandle $Q_{\tc}(L)$ is determined by two properties: the quandle relations in (\ref{defform}) and the fact that $Q_{\tc}(L)$ is translation-commutative. Using the notation $g_i = g_{x_i}$ for $1 \leq i \leq \mu$, the quandle relations in (\ref{defform}) provide the first kind of group relation mentioned in the statement of the corollary, and as noted in the proof of Proposition \ref{nilpotent}, the quandle relations in (\ref{tc}) provide group relations of the form $[g_i^{-1},[g_j^{-1},g_k^{-1}]]=1$. The $-1$ exponents in these commutators are not important, because elements of the form $[g_i^{-1},[g_j^{-1},g_k^{-1}]]$ generate the same normal subgroup of the free group on $\{g_1, \dots, g_ \mu \}$ as elements of the form $[g_i,[g_j,g_k]]$.
\end{proof}

For classical links with $\mu=2$, the first kind of relation in Corollary \ref{linknilpotent} was mentioned by Reidemeister in his famous monograph \cite[p.\ 45]{R}. Reidemeister's observation was extended to classical links of arbitrarily many components by Chen \cite{C}, who studied the lower central series of both $G(L)$ and the metabelian quotient $G(L)/G(L)''$; the lower central series quotients of $G(L)/G(L)''$ are called \emph{Chen groups} in his honor. Chen provided presentations of the groups $G(L)/G(L)_k$ for all $k \geq 3$. Milnor \cite{Mi} simplified Chen's presentations, and used them to define the $\bar \mu$-invariants. Only $G(L) /G(L) _3$ is important here, so this paper is not the right place to survey the extensive literature regarding these ideas. We refer to Hillman \cite{H} and Papadima and Suciu \cite{PS} for more general discussions, and further references to the classical theory. Recently, Chrisman \cite[Sec.\ 4]{Ch} has extended the Chen-Milnor theory of nilpotent quotients and $\bar \mu$-invariants to virtual links; although he does not single out $G(L) /G(L) _3$ for special attention, his results imply Corollary \ref{linknilpotent}.

Sakuma and Traldi \cite{ST} used Chen's presentation of $G(L)/G(L)_3$ to prove Theorem \ref{saktra} below. Recall that an \emph{articulation point} in a graph is a vertex whose removal results in a graph with strictly more connected components.

\begin{definition}
\label{insep}
Let $L=K_1 \cup \dots \cup K_{\mu}$ be a link, and let $\ell g (L)$ be the linking graph defined in Definition \ref{lg}. Then $L$ has \emph{inseparable linking numbers} if $\ell g (L)$ is connected, and has no articulation point. 
\end{definition}

For example, the link pictured in Fig.\ \ref{insepfig} of the introduction has inseparable linking numbers. No three-component sublink of this link has inseparable linking numbers, and every one- or two-component sublink has inseparable linking numbers. (In fact, every one- or two-component link has inseparable linking numbers.)

\begin{theorem}
\label{saktra}
(\cite{ST}) 
Let $L=K_1 \cup \dots \cup K_{\mu}$ and $L=K'_1 \cup \dots \cup K'_{\mu}$ be $\mu$-component classical links. Then these two properties are equivalent.
\begin{enumerate}
    \item There is an isomorphism between $G(L)/G(L)_3$ and $G(L')/G(L')_3$, which is ``meridian-preserving'' in the sense that for each $i \in \{1, \dots, \mu \}$, the image of a Wirtinger generator of $G(L)/G(L)_3$ corresponding to $K_i$ is a conjugate of a Wirtinger generator of $G(L')/G(L')_3$ corresponding to $K'_i$.
    \item Whenever $K_{i_1} \cup \dots \cup K_{i_ \nu} \subseteq L$ is a sublink with inseparable linking numbers, either $\ell(K_{i_j},K_{i_k})= \ell(K'_{i_j},K'_{i_k})$ $\thickspace \allowbreak \forall j \neq k \in \{1, \dots, \nu\}$ or $\ell(K_{i_j},K_{i_k})= -\ell(K'_{i_j},K'_{i_k})$ $\thickspace \allowbreak \forall j \neq k \in \{1, \dots, \nu \}$.
\end{enumerate}
\end{theorem}

The notion of ``meridian-preserving isomorphism'' used in Theorem \ref{saktra} can be described using quandles. There are natural maps $\eta:Q_{\tc}(L) \to \As(Q_{\tc}(L))$ and $\eta':Q_{\tc}(L') \to \As(Q_{\tc}(L'))$, given by $\eta(x)=g_x  \thickspace \allowbreak \forall x \in Q_{\tc}(L)$ and $\eta'(x')=g_{x'} \thickspace \allowbreak \forall x' \in Q_{\tc}(L')$. An isomorphism  $f:G(L)/G(L)_3 \to G(L')/G(L')_3$ is meridian-preserving if and only if for each $i \in \{1, \dots, \mu \}$,  $f(\eta(x_i))=\eta'(x')$ for some element $x'$ of the orbit of $x'_i$ in $Q_{\tc}(L')$.

Theorems \ref{thmtwo} and \ref{saktra} imply that $Q_{\tc}(L)$ (up to quandle isomorphisms) is a stronger invariant of classical links than  $G(L)/G(L)_3$ (up to re-indexings of $K_1, \dots,K_ \mu$ and meridian-preserving group isomorphisms). The reason is that Theorem \ref{thmtwo} involves entire connected components of $\ell g (L)$, while Theorem \ref{saktra} involves only connected subgraphs of $\ell g (L)$ that have no articulation points. These subgraphs can be much smaller than connected components.

\begin{figure} [b]
\centering
\begin{tikzpicture} [>=angle 90]
\draw [thick, domain=-30:180] plot ({-5+(0.8)*cos(\x)}, {(0.8)*sin(\x)});
\draw [thick, domain=180:295] plot ({-5+(0.8)*cos(\x)}, {(0.8)*sin(\x)});

\draw [thick, domain=150:180]  plot ({-3.95+(0.8)*cos(\x)}, {(0.8)*sin(\x)});
\draw [thick, domain=180:295] plot ({-3.95+(0.8)*cos(\x)}, {(0.8)*sin(\x)});

\draw [thick, domain=-30:115]  plot ({-3.95+(0.8)*cos(\x)}, {(0.8)*sin(\x)});

\draw [thick, domain=150:180] plot ({-2.9+(0.8)*cos(\x)}, {(0.8)*sin(\x)});
\draw [thick, domain=180:280] plot ({-2.9+(0.8)*cos(\x)}, {(0.8)*sin(\x)});
\draw [thick, domain=80:115] plot ({-2.9+(0.8)*cos(\x)}, {(0.8)*sin(\x)});


\draw [thick] [dashed] (-2.55,0) -- (-1.55, 0);

\draw [thick, domain=-30:80]  plot ({-1.3+(0.8)*cos(\x)}, {(0.8)*sin(\x)});
\draw [thick, domain=280:295] plot ({-1.3+(0.8)*cos(\x)}, {(0.8)*sin(\x)});

\draw [thick, domain=150:180]  plot ({-.25+(0.8)*cos(\x)}, {(0.8)*sin(\x)});
\draw [thick, domain=180:295] plot ({-.25+(0.8)*cos(\x)}, {(0.8)*sin(\x)});

\draw [thick, domain=-30:115]  plot ({-.25+(0.8)*cos(\x)}, {(0.8)*sin(\x)});
 \draw [thick, domain=150:180] plot ({.8+(0.8)*cos(\x)}, {(0.8)*sin(\x)});
\draw [thick, domain=180:475] plot ({.8+(0.8)*cos(\x)}, {(0.8)*sin(\x)});

\node at (-5,1.05) {$a_1$};
\node at (-3.95,1.05) {$a_2$};
\node at (-2.9,1.05) {$a_3$};
\node at (-.25,1.05) {$a_ {\mu-1}$};
\node at (-3.95,-1.05) {$b_2$};
\node at (-.25,-1.05) {$b_ {\mu-1}$};
\node at (-2.9,-1.05) {$b_3$};
\node at (.8,1.05) {$a_ \mu$};

\end{tikzpicture}
\caption{A connected sum of $\mu-1$ copies of the Hopf link.}
\label{cfig2}
\end{figure}
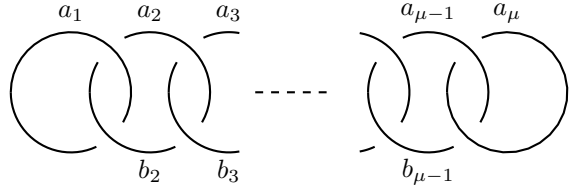

To illustrate this point, choose $\mu \geq 3$, and consider the oriented versions of the $\mu$-component link pictured in Fig.\ \ref{cfig2}. There are as many different oriented link types as there are lists of $\mu-1$ linking numbers all equal to $\pm 1$, with a list $\ell(K_1,K_2), \dots, \ell(K_{\mu-1},K_\mu)$ considered equivalent to the reversed list $\ell(K_{\mu-1},K_\mu), \dots, \ell(K_1,K_2)$. It follows that there are always at least 
\[
2+(2^{\mu-1}-2)/2= 2^{\mu-2}+1 
\]
oriented link types, because there are always at least two lists that remain the same when reversed, namely $1, \dots, 1$ and $-1, \dots, -1$. For each of these links, the entire graph $\ell g (L)$ is connected, but the only nontrivial connected subgraphs of $\ell g (L)$ without articulation points are single edges. Theorem \ref{thmtwo} tells us that the $Q_{\tc}$ quandles corresponding to two lists of linking numbers are isomorphic if and only if the components of the two links can be indexed so that the linking numbers differ only by a common factor of $ \pm 1$. Therefore, there are always at least
\[
\lceil (2^{\mu-2}+1)/2 \rceil =  2^{\mu-3}  +1 
\]
nonisomorphic $Q_{\tc}$ quandles among these links. In contrast, Theorem \ref{saktra} tells us that the nilpotent quotient groups $G(L)/G(L)_3$ are all isomorphic to each other, through meridian-preserving isomorphisms. In fact, even the $G(L)$ groups are all isomorphic through meridian-preserving isomorphisms:

\begin{proposition}
The groups of all oriented versions of the link pictured in Fig.\ \ref{cfig2} are isomorphic to each other, through meridian-preserving isomorphisms.
\end{proposition}

\begin{proof}
Let $L$ be one of these oriented links. Index the components of $L$ as $K_1, \dots, K_\mu$ in order from left to right, and index the arcs of the diagram as in Fig.\ \ref{cfig2}. For convenience, we use $g_a$ to denote the element of $G(L) = \As(Q(L))$ corresponding to an arc $a \in A(D)$, rather than $g_{q_a}$.

The two crossings involving $K_1$ and $K_2$ are on the far left. These two crossings provide two relations in $G(L)$; one relation is $g_{b_2}g_{a_1}g_{b_2}^{-1}= g_{a_1}$ and the other is either $g_{a_1}g_{a_2}g_{a_1}^{-1}= g_{b_2}$ or $g_{a_1}g_{b_2}g_{a_1}^{-1}= g_{a_2}$, depending on orientations. The first relation tells us that $g_{a_1}$ and $g_{b_2}$ commute with each other. Then either of the possible second relations tells us that $g_{a_2}=g_{b_2}$. With this equality in mind, one of the next two crossings tells us that $g_{b_3}g_{a_2}g_{b_3}^{-1}= g_{a_2}$, so $g_{a_2}$ and $g_{b_3}$ commute. The other crossing between $K_2$ and $K_3$ then provides a relation that guarantees $g_{a_3}=g_{b_3}$. With this equality in mind, the relations provided by the next two crossings tell us that $g_{a_3}$ and $g_{b_4}$ commute, and $g_{a_4}=g_{b_4}$. 

Continuing in this vein, we ultimately conclude that no matter how the link components are oriented, $G(L)$ is generated by $g_{a_1}, \dots, g_{a_{\mu}}$, subject to the relations $g_{a_i}g_{a_{i+1}}=g_{a_{i+1}}g_{a_i} \thickspace \allowbreak \forall i \in \{1, \dots, \mu-1 \}$. \end{proof}

\section{A theorem of Harrell and Nelson}
\label{hn}

Generally speaking, the absence of a structure theory for arbitrary quandles makes it difficult to work with them. One way to deal with this difficulty is to study special varieties of quandles, which can be described more easily. For instance, knot theorists have considered Alexander, involutory, latin, medial and $n$-quandles; see \cite{EN} for a survey. The translation-commutative quandles we have discussed in this paper constitute a subvariety of the medial quandles.

Another way to deal with the intractability of quandles is to consider numerical invariants derived from them. One way to derive such invariants of a link $L$ is to count the quandle maps from $Q(L)$ into specified target quandles. In \cite{HN}, Harrell and Nelson observed that if $L=K_1 \cup K_2$ is a two-component link, then there is a connection between the linking numbers $\ell_{1/2}(K_1,K_2), \ell_{2/1}(K_1,K_2)$ and the number of quandle maps from $Q(L)$ into quandles $X_n$ of a particular type. 

In the notation of Sec.\ \ref{tcq}, the quandle $X_n$ may be defined as follows. Suppose $n \geq 1 \in \mathbb Z$, let $B=\{b_1,b_2\}$, and let $S=\{S_{b_1},S_{b_2}\}$, where $S_{b_1}=\mathbb Z_B$ and $S_{b_2}$ is the subgroup of $\mathbb Z_B$ generated by $b_1$ and $nb_2$. Then $X_n=Q(S)$. Notice that $X_n$ is a translation-commutative quandle with a very simple structure: the first orbit has only one element, the second orbit has $n$ elements, $\beta_{b_1}$ cyclically permutes the elements of the second orbit, and $\beta_{b_2}$ is the identity map. Harrell and Nelson \cite{HN} proved the following.

\begin{proposition} (\cite{HN})
\label{HN}
Suppose $n \geq 2 \in \mathbb Z$. Let $L$ be a two-component link, and let $k$ be the number of quandle maps $Q(L) \to X_n$. 
\begin{enumerate}
    \item $k=n^2+1$ if $n$ divides neither $\ell_{1/2}(K_1,K_2)$ nor $\ell_{2/1}(K_1,K_2)$.
    \item $k=n^2+1+n$ if $n$ divides precisely one of $\ell_{1/2}(K_1,K_2), \ell_{2/1}(K_1,K_2)$.
    \item $k=n^2+1+2n$ if $n$ divides both $\ell_{1/2}(K_1,K_2)$ and $\ell_{2/1}(K_1,K_2)$.
\end{enumerate} 
\end{proposition}

As $X_n$ is translation-commutative, every quandle map $Q(L) \to X_n$ factors through the canonical surjection $Q(L) \to Q_{\tc}(L)$. Using this fact, it is not difficult to deduce Proposition \ref{HN} from Corollary \ref{qtc4}. For every two-component link $L$, there are $n^2+1$ non-surjective quandle maps $Q_{\tc}(L) \to X_n$, which are constant on each orbit of $Q_{\tc}(L)$. Moreover, if $i \neq j \in \{1,2\}$ and $n$ divides $\ell_{j/i}(K_1,K_2)$, then there are $n$ surjective quandle maps $Q_{\tc}(L) \to X_n$, which map $A_i = \mathbb Z_B / S_i(L)$ onto the $n$-element orbit of $X_n$ and map $A_j$ onto the one-element orbit of $X_n$.

\section*{Acknowledgments}
We are indebted to Kyle Miller for illuminating correspondence. We are also grateful to an anonymous reader, whose comments significantly improved the paper.


\begin{thebibliography}{99}

\bibitem{C} K.-T. Chen, Commutator calculus and link invariants, {\it Proc. Amer. Math. Soc.} {\bf 3} (1952) 44--55, 993.

\bibitem{Ch} M. Chrisman, Milnor's concordance invariants for knots on surfaces, to appear in {\it Algebr. Geom. Topol.}

\bibitem {EN} M. Elhamdadi and S. Nelson, {\it Quandles}, Student Mathematical Library, Vol. 74 (Amer. Math. Soc., Providence, R.I., 2015).

\bibitem {F} R. H. Fox, A quick trip through knot theory, in {\it Topology of 3-Manifolds and Related Topics} (Proc. The Univ. of Georgia Institute, 1961) (Prentice-Hall, Englewood Cliffs, N.J., 1962), pp. 120--167.

\bibitem{GPV} M. Goussarov, M. Polyak and O. Viro, Finite type invariants of classical and virtual knots, {\it Topology} {\bf 39} (2000) 1045--1068.

\bibitem {HN} N. Harrell and S. Nelson, Quandles and linking number, {\it J. Knot Theory Ramifications} {\bf 16} (2007) 1283--1293.

\bibitem {H} J. A. Hillman, {\it Algebraic Invariants of Links}, Series on Knots and Everything, Vol. 52, 2nd edn. (World Scientific, Singapore, 2012).

\bibitem {J} D. Joyce, A classifying invariant of knots, the knot quandle, {\it J. Pure Appl. Algebra} {\bf 23} (1982) 37--65.

\bibitem{K} G. Kuperberg, What is a virtual link?, {\it Algebr. Geom. Topol.} {\bf 3} (2003) 587--591.

\bibitem{MP} V. O. Manturov and D. P. Ilyutko, {\it Virtual Knots: The State of the Art}, Series on Knots and Everything, Vol. 51 (World Scientific, Singapore, 2013).

\bibitem {M} S. V. Matveev, Distributive groupoids in knot theory, {\it Mat. Sb. (N.S.)} {\bf 119} (1982) 78--88.

\bibitem {Mi} J. Milnor, Isotopy of links, in {\it Algebraic Geometry and Topology. A Symposium in Honor of S. Lefschetz}  (Princeton Univ. Press, Princeton, N.J., 1957), pp. 280--306.

\bibitem {PS} S. Papadima and A. Suciu, Chen Lie algebras, {\it Int. Math. Res. Not.} {\bf 2004} (2004) 1057--1086.

\bibitem {R} K. Reidemeister, {\it Knotentheorie} (Springer-Verlag, Berlin-New York, 1974).

\bibitem {RN} R. L. Ricca and B. Nipoti, Gauss' linking number revisited, {\it J. Knot Theory Ramifications} {\bf 20} (2011), 1325--1343.

\bibitem{ST} L. Traldi and M. Sakuma, Linking numbers and the groups of links, {\it Math. Semin. Notes, Kobe Univ.} {\bf 11} (1983) 119--132.

\end{thebibliography}
\end{document}